\documentclass[11pt,twoside]{article}
\usepackage{amsthm, amsmath, amssymb, amsfonts, graphicx, epsfig}

\usepackage{bbm}
\usepackage{mathrsfs}

\usepackage{epstopdf}

\usepackage{graphicx}
\usepackage[font=sl,labelfont=bf]{caption}
\usepackage{subcaption}

\usepackage{float}
\restylefloat{table}

\usepackage{enumerate}
\usepackage[colorlinks=true, pdfstartview=FitV, linkcolor=blue,
            citecolor=blue, urlcolor=blue]{hyperref}
\usepackage[usenames]{color}

\usepackage{url}
\makeatletter\def\url@leostyle{%
 \@ifundefined{selectfont}{\def\UrlFont{\sf}}{\def\UrlFont{\scriptsize\ttfamily}}} \makeatother

\usepackage[framemethod=default]{mdframed}

\usepackage[margin=1.0in, letterpaper]{geometry}

\newtheorem{theorem}{Theorem}

\newtheorem{lemma}[theorem]{Lemma}

\theoremstyle{definition}

\theoremstyle{remark}
\newtheorem{remark}[theorem]{Remark}

\numberwithin{equation}{section}
\numberwithin{theorem}{section}

\definecolor{Red}{rgb}{0.9,0,0.0}
\definecolor{Blue}{rgb}{0,0.0,1.0}

\def\cA{\mathcal{A}}

\def\cL{\mathcal{L}}

\def\cT{\mathcal{T}}

\def\bE{\mathbb{E}}
\def\bF{\mathbb{F}}

\def\bN{\mathbb{N}}

\def\bP{\mathbb{P}}
\def\bQ{\mathbb{Q}}
\def\bR{\mathbb{R}}

\def\sF{\mathscr{F}}

\newcommand{\Nz}{\mathbb N_0}

\newcommand{\set}[1]{\{#1\}}            
\renewcommand{\mid}{\;|\;}              

\DeclareMathOperator{\dif}{d \!}        

\title{ \vspace{-3em} 
   Risk-sensitive Markov decision problems under model uncertainty: finite time horizon case
}

\makeatletter
\def\and{%
  \end{tabular}%
  \begin{tabular}[t]{c}}%
\def\@fnsymbol#1{\ensuremath{\ifcase#1\or a\or b\or c\or
   d\or e\or f\or g\or h\or i\else\@ctrerr\fi}}
\makeatother

\author{
        Tomasz R. Bielecki\,\thanks{Department of Applied Mathematics, Illinois Institute of Technology
       \newline \hspace*{1.45em}  10 W 32nd Str, Building RE, Room 220, Chicago, IL 60616, USA
       \newline \hspace*{1.45em}  Emails: \url{tbielecki@iit.edu} (T. R. Bielecki), and \url{cialenco@iit.edu} (I. Cialenco)
       \newline \hspace*{1.45em}  URLs: \url{http://math.iit.edu/\~bielecki}  and \url{http://cialenco.com}
        \vspace{0.5em}} ,
\and
        Tao Chen\,\thanks{Department of Mathematics
University of Michigan
        \newline \hspace*{1.45em} 2074 East Hall,
530 Church Street,
Ann Arbor, MI 48109-1043, USA
        \newline \hspace*{1.45em} Email: \url{chenta@umich.edu}, URL: \url{http://taochen.im}
        \vspace{0.5em}},
\and
         Igor Cialenco,\,\footnotemark[1]        
        }

\date{ {\small This version: \today}} %

\begin{document}

\maketitle

{\footnotesize
\begin{tabular}{l@{} p{350pt}}
  \hline \\[-.2em]
  \textsc{Abstract}: \ & In this paper we study  a class of  risk-sensitive Markovian control problems in discrete time subject to model uncertainty. We consider a risk-sensitive discounted cost criterion with finite time horizon. The used methodology is the one of adaptive robust control combined with machine learning.
  \\[1em]
\textsc{Keywords:} \ & adaptive robust control, model uncertainty, stochastic control, adaptive robust dynamic programming, recursive confidence regions, risk-sensitive Markovian control problem, {machine learning, Gaussian surrogate processes, regression Monte Carlo.}
 \\[.5em]
\textsc{MSC2010:} \ & 49L20, 60J05, 60J20, 91A10, 91G10, 91G80, 62F25, 93C40, 93E35  \\[1em]
  \hline
\end{tabular}
}


\section{Introduction}
The main goal of this work is to study finite time horizon \textit{risk-sensitive Markovian control problems} subject to model uncertainty in a discrete time setup, and to develop a methodology to solve such problems efficiently. The proposed approach hinges on the following main building concepts: incorporating model uncertainty through the \textit{adaptive robust} paradigm introduced in \cite{BCCCJ2019} and developing efficient numerical solutions for the obtained Bellman equations by adopting the \textit{machine learning techniques} proposed in \cite{ChenLudkovski2019}.

There exists a significant body of work on incorporating model uncertainty (or model misspecification) in stochastic control problems, and among some of the well-known and prominent methods we would mention the robust control approach \cite{GilboaSchmeidler1989,HansenSargent2006,HansenBookBook2008}, adaptive control \cite{KumarVaraiya2015Book,ChenGuo1991-Book}, and Bayesian adaptive control \cite{KumarVaraiya2015Book}. A comprehensive literature review on this subject is beyond the scope of this paper, and we refer the reader to \cite{BCCCJ2019} and references therein. In~\cite{BCCCJ2019} the authors proposed a novel adaptive robust methodology that solves time-consistent Markovian control problems in discrete time subject to model uncertainty - the approach that we take in this study too. The core of this methodology was to combine a recursive learning mechanism about the unknown model with the underlying Markovian dynamics, and to demonstrate that the so called adaptive
robust Bellman equations produce an optimal adaptive robust control strategy.

In contrast to \cite{BCCCJ2019}, where  the considered optimization criterion was of the terminal reward type, in the present work, we also allow intermediate rewards and we use the discounted risk sensitive criterion. Accordingly, we derive a new set of adaptive robust Bellman equations, similar to those used in \cite{BCCCJ2019}.

Risk sensitive criterion has been broadly used both in the control oriented  literature, as well as in the game oriented literature. We refer to, e.g.,  \cite{Bielecki2003,DavisLleo2014,BauRie2017}, and the references  therein for insight into risk sensitive control and risk sensitive games both in discrete time and in continuous time.

The paper is organized as follows. In Section~\ref{sec:robust} we formulate the finite time horizon risk-sensitive Markovian control problem subject to model uncertainty that is studied here. Section~\ref{sec:formulation} is devoted to the formulation and to study of the robust adaptive control problem that is relevant for the problem formulated in Section~\ref{sec:robust}. This section  presents the main theoretical developments of the present work. 
In Section~\ref{sec:ex} we formulate an illustrative example of our theoretical results that is rooted in the classical linear-quadratic-exponential control problem (see e.g. \cite{HanSar1995}). Next, using machine learning methods, in Section~\ref{sec:ml} we provide numerical solutions of the example presented in Section~\ref{sec:ex}.


Finally, we want to mention that the important case of an infinite time horizon risk-sensitive Markovian control problem in discrete time subject to model uncertainty will be studies in a follow-up work.

\section{Risk-sensitive Markovian discounted control problems with model uncertainty}\label{sec:robust}
In this section we state the  underlying discounted risk-sensitive stochastic control problems.
Let  $(\Omega, \sF)$ be a measurable space,  $T\in \bN$ be a finite time horizon, and let us denote by $\cT:=\set{0,1,2,\ldots,T}$ and $\cT':=\set{0,1,2,\ldots,T-1}$. We let $\boldsymbol \Theta\subset \bR^d$   be a non-empty compact set, which will play the role of the parameter space throughout. We consider a random process $Z=\{Z_t,\ t=1,2\ldots\}$ on $(\Omega, \sF)$  taking values in $\bR^m$, and we denote by ${\mathbb {F}}=(\sF_t,\ t=0,2\ldots)$ its natural filtration, with $\sF_0=\set{\emptyset,\Omega}$. We postulate that this process is observed by the controller, but the true law of $Z$  is unknown to the controller and assumed to be generated by a  probability measure belonging to a (known) parameterized family of probability distributions on $(\Omega, \sF)$, say $\mathbf{P}(\boldsymbol \Theta)=\{\bP_\theta,  \theta\in \boldsymbol \Theta\}$.  As usually, $\bE_\bP$ will denote the expectation under a probability measure $\bP$ on  $(\Omega, \sF)$, and, for simplicity, we will write $\bE_\theta$ instead of  $\bE_{\bP_\theta}$. We denote by $\bP_{\theta^*}$ the measure generating the true law of $Z$, and thus $\theta ^*\in \boldsymbol
\Theta$ is the unknown true parameter. The sets $\boldsymbol \Theta$ and  $\mathbf{P}(\boldsymbol \Theta)$ are known to the observer. Clearly, the model uncertainty may occur if $\boldsymbol \Theta \ne \{\theta^*\}$, which we will assume to hold throughout.

We let $A\subset \bR^k$ be a finite set,\footnote{$A$ will represent the set of control values, and we assume it is finite for simplicity, in order to avoid technical issues regarding the existence of measurable selectors.}   and $S\, :\, \bR^n\times A\times \bR^m \rightarrow \bR^d$ be a measurable mapping.
An admissible control process $\varphi$ is an $\bF$-adapted process, taking values in $A$, and we will denote by $\cA$ the set of all admissible control processes.

We consider an  underlying discrete time controlled dynamical system with the state process $X$ taking values in $\bR^n$ and control process $\varphi$ taking values in $A$. Specifically, we let
\begin{equation}\label{eq:mm}
X_{t+1}=S(X_t,\varphi_t,Z_{t+1}),\quad t\in \cT' ,\quad X_0=x_0 \in \bR^n.
\end{equation}

At each time $t=0,\ldots,T-1,$ the running reward $r_t(X_t,\varphi_t)$ is delivered, where, for every $a\in A$, the function $r_t(\cdot,a): \bR^n\to\bR_+$ is bounded and continuous. Similarly, at the terminal time $t=T$ the terminal reward $r_T(X_T)$ is delivered, where $r_T: \bR^n\to\bR_+$ is a bounded and continuous function.

Let $\beta \in (0,1)$ be a discount factor, and let $\gamma \ne 0$ be the risk sensitivity factor. The underlying discounted, risk-sensitive control problem is:
\begin{equation}\label{eq:problem-T}
\sup_{\varphi\in \cA} \frac{1}{\gamma}\ln\left (\bE_{\theta^*} e^{\gamma \left(\sum_{t=0}^{T-1}\beta^t r_t(X_t,\varphi_t)+\beta^T r_T(X_T)\right)}\right )
\end{equation}
subject to \eqref{eq:mm}. Clearly, since $\theta^*$ is not known to the controller, the above problem can not be solved as it is stated. The main goal of this paper is formulate and solve the adaptive robust control problem corresponding to \eqref{eq:problem-T}.

\begin{remark}
(i) The risk-sensitive criterion in \eqref{eq:problem-T} is in fact an example of application of the entropic risk measure, say $\rho_{\theta^*,\gamma}$, which is defined as
\[
\rho_{\theta^*,\gamma}(\xi):= \frac{1}{\gamma}\ln\bE_{\theta^*} e^{\gamma \xi},
\]
where $\xi$ is a random variable on $(\Omega, \sF,P_{\theta^*})$ that admits finite moments of all orders.
\\
(ii) It can be verified that
\begin{align*}
\rho_{\theta^*,\gamma}(\xi)=\bE_{\theta^*}(\xi) + \frac {\gamma}{2}{\mathbb{VAR}}_{\theta^*}(\xi) +O(\gamma^2).
\end{align*}
Thus, in case when $\gamma < 0$ the term $\frac {\gamma}{2}{\mathbb{VAR}}_{\theta^*}(\xi)$ can be interpreted as the risk-penalizing term. On the contrary, when $\gamma > 0$, the term $\frac {\gamma}{2}{\mathbb{VAR}}_{\theta^*}(\xi)$ can be viewed as the risk-favoring term.\\
(iii) In the rest of the paper we focus on the case $\gamma >0$. The case $\gamma<0$ can be treated in an analogous way.
\end{remark}

\section{The adaptive robust risk sensitive  discounted control problem}\label{sec:formulation}

We follow here the developments presented in \cite{BCCCJ2019}. The key difference is that in this work we deal with running and terminal costs.

In what follows, we will be making use of a recursive construction of confidence regions for the unknown parameter $\theta^*$ in our model. We refer to \cite{BCC2017} for a general study of recursive constructions of (approximate) confidence regions for time homogeneous Markov chains. Section~\ref{sec:ex} provides details of a specific  such recursive construction corresponding to the example presented in that section.   Here, we just postulate that the recursive algorithm for building confidence regions uses a $\boldsymbol{\Theta}$-valued and observed  process, say $C=(C_t,\ t\in \Nz)$, satisfying the following abstract dynamics
\begin{equation}\label{eq:R}
C_{t+1}= R(t,C_t,Z_{t+1}),\quad t\in \Nz,\ C_0=c_0\in\boldsymbol{\Theta},
\end{equation}
where $R:\Nz\times\bR^d\times\bR^m \to \boldsymbol{\Theta}$ is a deterministic measurable function.
Note that, given our assumptions about process $Z$, the process $C$ is $\bF$-adapted.
This is   one of the key features of our model. Usually $C_t$ is taken to be a consistent estimator of $\theta^*$.

Now, we fix a confidence level $\alpha\in (0,1),$ and for each time $t\in \Nz$, we assume that  an $(1-\alpha)$-confidence region, say $\mathbf{\Theta}_t \subset \bR^d$, for $\theta^*$, can be represented as
\begin{equation}\label{eq:CIR}
\mathbf{\Theta}_t=\tau(t,C_t),
\end{equation}
where, for each $t\in \Nz$, $\tau(t,\cdot)\, :\, \bR^d \rightarrow 2^{\mathbf{\Theta}}$ is a deterministic set valued   function, where, as usual, $2^{\mathbf{\Theta}}$ denotes the set of all subsets of ${\mathbf{\Theta}}$. Note that in view of \eqref{eq:R} the construction of confidence regions given in \eqref{eq:CIR} is indeed recursive. In our construction of confidence regions, the mapping $\tau(t,\cdot)$ will be a measurable set valued function,  with compact values. It needs to be noted that we will only need to compute $\mathbf{\Theta}_t$ until time $T-1$.
In addition, we assume that for any $t\in\cT'$, the mapping $\tau(t,\cdot)$ is upper hemi-continuous (u.h.c.). That is, for any $c\in\boldsymbol\Theta$, and any open set $E$ such that $\tau(t,c)\subset E \subset\boldsymbol\Theta$, there exists a neighbourhood $D$ of $c$ such that for all $c'\in D$, $\tau(t,c')\subset E$ (cf. \cite[Definition~11.3]{Border1985}).

\begin{remark}
 The important property of the recursive confidence regions constructed as indicated above is that, in many models, $\lim_{t\rightarrow \infty} \mathbf{\Theta}_t=\set{\theta^*}$, where the convergence is understood $\bP_{\theta^*}$ almost surely, and the limit is in the Hausdorff metric. This is not always the case though in general. In \cite{BCC2017} is shown that the convergence holds in probability, for the model setup studied there.
\end{remark}

The sequence $\mathbf{\Theta}_t,\ t\in \cT'$ represents learning about $\theta^*$ based on the observation of the history $(Y_0,Y_1\ldots,Y_t),\ t\in \cT'$, where  $Y_t=(X_t,C_t),\ t\in\cT,$ is the  augmented state process taking values in the augmented state space
\[
E_Y =  \bR^n  \times \boldsymbol{\Theta}.
\]
We denote by ${\mathcal E}_Y$ the collection of Borel measurable sets in $E_Y$.

In view of the above, if the control process $\varphi$ is employed then the process $Y$ has the following dynamics
\[
Y_{t+1}=\mathbf{G}(t, Y_t,\varphi_t,Z_{t+1}),\quad t\in \cT',
\]
where  the mapping  $\mathbf{G}\, :\, \Nz\times E_Y \times A \times \bR^m \rightarrow E_Y$ is defined as
\begin{equation}\label{eq:T}
\mathbf{G}(t,y,a,z)=\big(S(x,a,z), R(t,c,z)\big),
\end{equation}
with $ y=(x,c)\in  E_Y$.

We define the corresponding histories
\begin{equation}\label{hist}
H_t=(Y_0,\ldots,Y_t),\quad t\in \cT',
\end{equation}
so that
\begin{equation}\label{eq:boldh}
H_t\in \mathbf{H}_t= \underbrace{E_Y \times E_Y \times \ldots
\times E_Y}_{t+1 \textrm{ times}}.
\end{equation}
Clearly, for any admissible control process $\varphi$, the random variable $H_t$ is $\sF_t$-measurable. We denote by
\begin{equation}\label{eq:h}
h_t=(y_0,y_1,\ldots,y_t)=(x_0,c_0,x_1,c_1,\ldots,x_t,c_t)
\end{equation}
 a realization of $H_t.$ Note that $h_0=y_0=(x_0,c_0)$.

A control process $\varphi=(\varphi_t,\ t \in \cT')$ is called  history dependent control process if (with a slight abuse of notation)
\[
\varphi_t = \varphi_t(H_t) ,
\]
where (on the right hand side) $\varphi_t \, :\, \mathbf{H}_t \rightarrow A$, is a measurable mapping. Given our above setup, any history dependent control process is ${\mathbb {F}}$--adapted, and thus, it is admissible. For any admissible control process $\varphi$ and for any $t\in\cT'$, we denote by $\varphi^{t}=(\varphi_k,\ k=t,\dots, T-1)$ the `$t$-tail' of $\varphi$. Accordingly, we denote by ${\mathcal A}^{t}$  the collection of `$t$-tails' of $\varphi$. In particular, $\varphi^{0}=\varphi$ and ${\mathcal A}^{0}={\mathcal A}$. The superscript notation applied to processes should not be confused with power function applied such as $\beta^t$.

Let $\psi_t:\mathbf{H}_t\to\mathbf{\Theta}$  be a Borel measurable mapping such that $\psi_t(h_t)\in \tau(t,c_t)$, and let us denote by $\psi =(\psi_t, \ t\in\cT')$ the sequence of such mappings, and by $\psi^{t}$  the $t$-tails of the sequence $\psi$, in analogy to $\varphi^{t}$. The set of all sequences $\psi$, and respectively $\psi^{t}$ , will be denoted by $\mathbf{\Psi}$ and $\mathbf{\Psi}^{t}$, respectively.

Strategies $\varphi$ and $\psi$ are called \textit{Markovian strategies or policies} if (with some abuse of notation)
\[
\varphi_t = \varphi_t(Y_t) ,\quad \psi_t = \psi_t(Y_t),
\]
where (on the right hand side) $\varphi_t \, :\, E_Y \rightarrow A$, and  is a (Borel) measurable  mapping, and  $\psi_t \, :\, E_Y \rightarrow \mathbf{\Theta}$ is a (Borel) measurable mapping satisfying  $\psi_t(x,c)\in \tau(t,c)$.

 In order to simplify all the following argument we limit ourselves to Markovian policies. In case of Markovian dynamics settings, such as ours, this comes without loss of generality, as there typically exist optimal Markovian strategies, if optimal strategies exist at all. Accordingly, $\mathcal A$ and $\mathbf{\Psi}$ are now sets of Markov strategies.

Next, for each $(t,y,a,\theta)\in  \cT'\times  E_Y\times A\times \mathbf{\Theta}$,  we define a probability measure on $\mathcal{E}_Y$, given by
\begin{equation}\label{eq:QB-bar}
 Q(B\mid t,y,a,\theta)=\bP_\theta(Z_{t+1}\in \{z: \mathbf{G}(t,y,a,z)\in B\})=\bP_\theta\left(\mathbf{G}(t,y,a,Z_{t+1})\in B\right),\ B\in \mathcal{E}_Y.
\end{equation}
We assume that for every $t\in\cT$ and every $a\in A$,  we have that $Q(dy'\mid t,y,a,\theta)$  is a Borel measurable stochastic kernel with respect to $(y,\theta)$. This assumption will be strengthened later on.

Using Ionescu-Tulcea theorem (cf. \cite[Appendix B]{Baeuerle2011book}), for every $t=0,\ldots,T-1$, every $t$-tail $\varphi^{t} \in \cA^{t} $ and every state $y_t\in E_Y$,  we define the family ${\mathcal Q}^{\varphi^{t},{\boldsymbol \Psi^{t}}}_{y_t,t} =\{{\mathbb Q}^{\varphi^{t},\psi^{t}}_{y_t, t},\  \boldsymbol\psi^{t} \in {\boldsymbol \Psi^{t}} \}$  of probability measures on the  concatenated canonical space $\textsf{X}_{s=t+1}^T E_Y$, with
\begin{align}\label{eq:prob}
\mathbb{Q}^{\varphi^{t}, {\psi}^{t}}_{y_t,t}(B_{t+1}\times \cdots \times B_T) :=
\int\limits_{B_{t+1}}\cdots \int\limits_{B_T}
\prod\limits_{u=t+1}^{T} Q( \dif y_u\mid u-1,y_{u-1},\varphi_{u-1}(y_{u-1}),{\psi}_{u-1}(y_{u-1})).
\end{align}

The \textit{discounted, risk-sensitive, adaptive robust control problem} corresponding\footnote{Since $\gamma>0$, we omit the factor $1/\gamma$.} to \eqref{eq:problem-T} is:
\begin{equation}\label{eq:problem-T-robust-noln}
\sup_{\varphi^{0}\in \cA^{0}} \inf_{\mathbb{Q}\in {\mathcal Q}^{\varphi^{0},{\boldsymbol \Psi}^{0}}_{y_0,0}}
\bE_{\mathbb Q} e^{\gamma\sum_{t=0}^{T}\beta^t r_t(X_t,\varphi_t(Y_t))},
\end{equation}
where, for simplicity of writing,  here and everywhere below, with slight abuse of notations,  we set $r_T(x,a)=r_T(x)$.
In next section we will show that a solution to this problem can be given in terms of the discounted adaptive robust Bellman equations associated to it.

\subsection{Adaptive robust Bellman equation}\label{sec:adaptRobustBellmanEq}

Towards this end we aim our attention at the following adaptive robust Bellman equations
\begin{align}
  W_{T}(y) & = e^{\gamma \beta^T r_T(x)}, \quad y\in E_Y,
    \nonumber \\
W_{t}(y) & = \max_{a\in A} \inf_{ \theta \in \tau(t,c)}\int_{E_Y} W_{{t+1}}(y' )e^{\gamma \beta^t r_t(x,a)}
   Q( \dif y'\mid t,y,a,\theta ) , \quad
     y\in E_Y, \   t=T-1, \ldots, 0, \label{eq:bellmanEquationrobustIII}
\end{align}
where we recall that $y=(x,c)$.
\begin{remark}
Clearly, in \eqref{eq:bellmanEquationrobustIII}, the exponent $e^{\gamma \beta^t r_t(x,a)}$ can be factored out, and $W_t$ can be written as
$$
W_{t}(y) = \max_{a\in A} \left( e^{\gamma \beta^t r_t(x,a)} \cdot \inf_{ \theta \in \tau(t,c)}\int_{E_Y} W_{{t+1}}(y' )
Q( \dif y'\mid t,y,a,\theta ) \right).
$$
Nevertheless, in what follows, we will keep similar factors inside of the integrals, mostly for the convenience of writing as well as to match the visual appearance of classical Bellman equations.
\end{remark}

We will study the solvability of this system. We start with Lemma~\ref{lemma:WTusc} below, where, under some additional technical assumptions, we show that the optimal selectors in \eqref{eq:bellmanEquationrobustIII} exist; namely, for any $t\in\cT'$,  and any $y=(x,c)\in E_Y$, there exists a measurable mapping $\varphi^*_t\, : \, E_Y \rightarrow A$, such that
\[
 W_{t}(y)  =  \inf_{ \theta \in \tau(t,c)}\int_{E_Y}W_{{t+1}}(y' )e^{\gamma \beta^t r_t(x, \varphi_t^*(y) )}
   Q( \dif y'\mid t,y,\varphi^{*}_t(y),\theta).
\]
In order to proceed, for the sake of simplicity, we will assume that under measure $\bP_\theta$, for each $t\in \cT$, the random variable $Z_t$ has a density   with respect to the Lebesgue measure, say $f_Z(z;\theta), \ z\in\bR^m$. In this case we have
\[
\int_{E_Y}W_{{t+1}}(y' )
   Q( \dif y'\mid t,y,a,\theta)=\int_{\bR^m} W_{{t+1}}(\mathbf{G}(t,y, a, z))f_Z(z; \theta)\dif z,
\]
where $\mathbf{G}(t,y, a, z)$ is given in \eqref{eq:T}.

Additionally, we take the standing assumptions:
\begin{enumerate}[(i)]
\item for any $a$ and $z$, the function $S(\cdot,a,z)$ is
continuous;
  \item for each $z$, the function $f_Z(z;\cdot)$ is continuous; 
  \item for each $t\in \cT'$, the function $R(t,\cdot,\cdot)$ is continuous.
\end{enumerate}
Then, the following result holds true.

\begin{lemma}\label{lemma:WTusc} The functions $W_{t},\ t=T,T-1,\ldots,0,$ are lower semi-continuous (l.s.c.), and the optimal selectors $\varphi^{*}_t ,\ t=T-1,\ldots,0,$ realizing maxima in \eqref{eq:bellmanEquationrobustIII} exist.
\end{lemma}

\begin{proof} Since $r_T$ is continuous and bounded, so is the function $W_{T}$.
Since $\mathbf{G}(T-1,\cdot,a,z)$ is continuous, then, $W_{T}(\mathbf{G}(T-1,\cdot,a,z))$ is continuous.
Consequently, recalling again that $y=(x,c)$, for each $a$, the function
\[
w_{T-1}(y,a,\theta) := \int_{\bR} W_{T}(\mathbf{G}(T-1,y, a, z))e^{\gamma \beta^{T-1} r_{T-1}(x,a)} f_Z(z; \theta)\dif z
\]
\[=e^{\gamma \beta^{T-1} r_{T-1}(x,a)}\int_{\bR} e^{\gamma \beta^{T} r_T(S(x,a,z))}f_Z(z; \theta)\dif z
\]
is continuous in $(y,\theta)$.

Next, we will apply \cite[Proposition 7.33]{BertsekasShreve1978Book} by taking (in the notations of \cite{BertsekasShreve1978Book})
\begin{align*}
\mathrm{X} & = E_Y\times A = \bR^n\times \boldsymbol{\Theta} \times A,  \quad \mathrm{x}=(y,a), \\
\mathrm{Y} & =  \boldsymbol{\Theta}, \quad  \mathrm{y}=\theta, \\
\mathrm{D} & =\bigcup_{(y,a)\in E_Y \times A} \set{(y,a)} \times  \tau (T-1,c), \\
f(\mathrm{x},\mathrm{y}) & =  w_{T-1}(y,a,\theta).
\end{align*}
Note that in view of the prior assumptions, $\mathrm{Y}$ is metrizable and compact.
Clearly $\mathrm{X}$ is metrizable.
From the above, $f$ is continuous, and thus lower semi-continuous.  Since $\tau(T-1,\cdot)$ is compact-valued and u.h.c. on $E_Y\times A$, then according to \cite[Proposition 11.9]{Border1985}, the set-valued function $\tau(T-1,\cdot)$ is closed, which implies that its graph $\mathrm{D}$ is closed \cite[Definition~11.5]{Border1985}. Also note that the cross section $\mathrm{D}_\mathrm{x} = {\mathrm D}_{(y,a)} = \set{\theta \in \boldsymbol{\Theta}\, :\, (y,a,\theta) \in \mathrm{D}}$ is given by $\mathrm{D}_{(y,a)}(T-1)=\tau (T-1,c)$.
Hence, by \cite[Proposition 7.33]{BertsekasShreve1978Book}, the function
\[
\widetilde w_{T-1}(y,a)= \inf_{\theta \in \tau(T-1,c)} (w_{T-1}(y,a,\theta)) ,\quad (y,a)\in E_Y \times A,
\]
is l.s.c..  Consequently, the function $\widehat w_{T-1}(y,a)=-\widetilde w_{T-1}(y,a)$ is upper semi-continuous (u.s.c). Thus, by \cite[Proposition 7.34]{BertsekasShreve1978Book}, the function  $-W_{T-1}(y)= -{\max _{a\in A}}\widetilde w_{T-1}(y,a)={\min _{a\in A}}\widehat w_{T-1}(y,a)$ is u.s.c., so that $W_{T-1}(y)$ is l.s.c.. Moreover, since $A$ is finite, there exists an optimal selector $\varphi^*_{T-1}$, that is $W_{T-1}(y)=\widetilde w_{T-1}(y,\varphi^*_{T-1}(y))$.

Proceeding to the next step, note that $W_{T-1}(\mathbf{G}(T-2,y,a,z))e^{\gamma\beta^{T-2}r_{T_2}(x,a)}$ is l.s.c. and positive, hence bounded from below.
Therefore, according to \cite[Proposition 7.31]{BertsekasShreve1978Book}, the function
$$
w_{T-2}(y,a,\theta)=\int_\bR W_{T-1}(\mathbf{G}(T-2,y,a,z))e^{\gamma\beta^{T-2}r_{T-2}(x,a)}f_Z(z;\theta)\dif z
$$
is l.s.c..
The rest of the proof follows in the analogous way.
\end{proof}


Next, we will prove an auxiliary result needed to justify the mathematical operations conducted in the proof of the main result -- Theorem~\ref{th:main}.
Define the functions $U_{t}$ and $U^*_{t}$ as follows: for $\varphi^{t}\in {\mathcal A}^{t}$ and $y  \in E_Y$,
\begin{align}\label{eq:prob1R-adaptive-Dirac}
   U_{t}(\varphi^{t},y)
   &= e^{\gamma \beta^t r_t(x,\varphi_t(y))}\inf_{\mathbb{Q}\in {\mathcal Q}^{\varphi^{t},{\boldsymbol \Psi}^{t}}_{y,t}}
   \bE_{\mathbb Q} e^{\gamma \sum_{k=t+1}^{T}\beta^k r_k(X_k,\varphi_k(Y_k))}, \quad t\in\cT',\\
    U^*_{t}(y)&=\sup_{\varphi^{t}\in {\mathcal A}^{t}} U_{t}(\varphi^{t},y),\quad t\in\cT',\\
    U^*_{T}(y) & = e^{\gamma \beta ^T r_T(x)}.
\end{align}

We now have the following result.

\begin{lemma}\label{lemma:U-lsa} For any $t\in \cT'$, and for any $\varphi^{t}\in {\mathcal A}^{t}$, the function $U_{t}(\varphi^{t},\cdot)$ is lower semi-ananlytic (l.s.a.) on $E_Y$. Moreover, there exists a sequence of universally measurable functions $\psi^*_k$, $k=t,\ldots,T-1$ such that
\begin{align}\label{eq:lemma34}
U_t(\varphi^t,y)=e^{\gamma\beta^tr_t(x,\varphi_t(y))}\bE_{\bQ^{\varphi^t,\psi^{t,*}}_{y,t}}e^{\gamma\sum_{k=t+1}^T\beta^kr_k(X_k,\varphi_k(Y_k))}.
\end{align}
\end{lemma}

\begin{proof}  According to \eqref{eq:QB-bar}, and using the definition of ${\mathcal Q}^{\varphi^{t},{\boldsymbol \Psi}^{t}}_{y,t}$,  we have that
\begin{align}
U_{t}(\varphi^{t},y) & =\inf_{\psi^{t}\in\boldsymbol\Psi^{t}}\int_{E_Y}\cdots\int_{E_Y}
e^{\gamma\sum_{k=t}^{T}\beta^kr_k(x_k,\varphi_k(y_k))}
Q(\dif y_T|T-1,y_{T-1},\varphi_{T-1}(y_{T-1}),\psi_{T-1}(y_{T-1})) \nonumber \\
& \qquad \qquad \qquad \qquad \cdots Q( \dif y_{t+1}|t,y,\varphi_{t}(y),\psi_{t}(y)). \label{eq:U3}
\end{align}
For a given policy $\varphi\in\cA$, define the following functions on $E_Y$
\begin{align*}
V_{T}(y)&=e^{\gamma \beta^T r_T(x)},\\
V_{t}(y)&=\inf_{\theta\in\tau(t,c)}\int_{E_Y}e^{\gamma\beta^tr_t(x,\varphi_t(y)}V_{t+1}(y')Q(\dif y'|t,y,\varphi_t(y),\theta),\quad t\in\cT'.
\end{align*}
We will prove recursively that the functions $V_{t}$ are l.s.a. in $y$, and that
\begin{equation}\label{V}
V_{t}(y)=U_{t}(\varphi^{t},y),\quad t=0,\ldots,T-1.
\end{equation}
Clearly, $V_{T}$ is l.s.a. in $y$.

Next, we will prove that $V_{T-1}(y)$ is l.s.a.. By our assumptions, the stochastic kernel $Q(\cdot|T-1,\cdot,\cdot,\cdot)$ is Borel measurable on $E_Y$ given $E_Y\times A\times \boldsymbol\Theta$, in the sense of \cite[Definition~7.2]{BertsekasShreve1978Book}.
Then, the integral $\int_{E_Y}V_T(y')  Q(\dif y'|T-1,y,a,\theta)$ is l.s.a. on $E_Y\times A\times\boldsymbol\Theta$ according to \cite[Proposition 7.48]{BertsekasShreve1978Book}.
Now, we set (in the notations of \cite{BertsekasShreve1978Book})
\begin{align*}
\mathrm{X}&=E_Y\times A,\quad \mathrm{x}=(y,a)\\
\mathrm{Y}&=\boldsymbol\Theta, \quad \mathrm{y}=\theta,\\
\mathrm{D}&=\bigcup_{(y,a)\in E_Y\times A}\{y,a\}\times\tau(T-1,c),\\
f(\mathrm{x},\mathrm{y})&=\int_{E_Y}V_T(y') \, Q(\dif y'|T-1,y,a,\theta).
\end{align*}
Note that in view of our assumptions, $\mathrm{X}$ and $\mathrm{Y}$ are Borel spaces.
The set $\mathrm{D}$ is closed (see the proof of Lemma~\ref{lemma:WTusc}) and thus analytic.
Moreover, $\mathrm{D}_{\mathrm{x}}=\tau(T-1,c)$.
Hence, by \cite[Proposition 7.47]{BertsekasShreve1978Book}, for each $a\in A$ the function
$$
\inf_{\theta\in\tau(T-1,c)}\int_{E_Y} V_T(y') \, Q(\dif y'|T-1,y,a,\theta)
$$
is l.s.a. in $y$. Thus, it is l.s.a. in $(y,a)$.
Moreover, in view of \cite[Proposition~7.50]{BertsekasShreve1978Book}, for any $\varepsilon>0$, there exists an analytically measurable function $\psi^{\varepsilon}_{T-1}(y,a)$ such that
$$
\inf_{\theta\in\tau(T-1,c)}\int_{E_Y} V_T(y') \, Q(\dif y'|T-1,y,a,\theta)=\int_{E_Y} V_T(y') \, Q(\dif y'|T-1,y,a,\psi^{\varepsilon}_{T-1}(y,a))+\varepsilon.
$$
Therefore, for any fixed $(y,a)$, we obtain a sequence $\{\psi^{1/n}_{T-1}(y,a),n\in\bN\}$ such that
$$
\lim_{n\to\infty}\int_{E_Y} V_T(y') \, Q(\dif y'|T-1,y,a,\psi^{1/n}_{T-1}(y,a))=\inf_{\theta\in\tau(T-1,c)}\int_{E_Y} V_T(y') \, Q(\dif y'|T-1,y,a,\theta).
$$
Due to the assumption that $\tau(T-1,c)$ is compact, there exists a convergent subsequence $\{\psi^{1/n_k}_{T-1}(y,a),k\in\bN\}$ such that its limit $\psi^*_{T-1}(y,a)$ is universally measurable and satisfies
$$
\int_{E_Y} V_T(y') \, Q(\dif y'|T-1,y,a,\psi^*_{T-1}(y,a))=\inf_{\theta\in\tau(T-1,c)}\int_{E_Y} V_T(y') \, Q(\dif y'|T-1,y,a,\theta).
$$
Clearly, the function $e^{\gamma\beta^{T-1}r_{T-1}(x,a)}$ is l.s.a. in $(y,a)$.
Thus, since  $\varphi_{T-1}(y)$ is a Borel measurable function, using part (3) in \cite[Lemma~7.30]{BertsekasShreve1978Book} we conclude that
both $e^{\gamma\beta^{T-1}r_{T-1}(x,\varphi_{T-1}(y))}$ and $\inf_{\theta\in\tau(T-1,c)}\int_{E_Y} V_T(y')Q(\dif y'|T-1,y,\varphi_{T-1}(y),\theta)$ are l.s.a. in $y$. Since both these functions are non-negative then, by part (4) in \cite[Lemma 7.30]{BertsekasShreve1978Book}, we conclude that $V_{T-1}$ is l.s.a. in $y$. The proof that $V_{t}$ is l.s.a. in $y$ and $\psi^*_t$ exists for $t=0,\ldots, T-2$, follows analogously. We also obtain that
\begin{align}\label{eq:lemma34-3}
\int_{E_Y} V_t(y') \, Q(\dif y'|t-1,y,a,\psi^*_{t-1}(y,a))=\inf_{\theta\in\tau(t-1,c)}\int_{E_Y} V_t(y') \, Q(\dif y'|t-1,y,a,\theta),
\end{align}
for any $t=1,\ldots,T-1$.

It remains to verify \eqref{V}.  For $t=T-1$, by \eqref{eq:U3}, we have
\begin{align*}
U_{T-1}(\varphi^{T-1},y)
& =\inf_{\theta\in\tau(T-1,c)}\int_{E_Y}e^{\gamma\beta^{T-1}r_{T-1}(x,\varphi_{T-1}(y))}V_{T}(y')Q(\dif y'|T-1,y,\varphi_{T-1}(y),\theta)\\
& =V_{T-1}(y).
\end{align*}
Therefore, $U_{T-1}(\varphi^{T-1},\cdot)$ is l.s.a..
Assume that for $t=1,\ldots,T-1$, $U_{t}(\varphi^{t},y)=V_{t}(y)$, and it is l.s.a..
Then, for any $y_{t-1}\in E_Y$, with the notation $\psi^{t-1} = (\psi_{t-1}, \psi^t)$,  we get
\begin{align*}
&U_{t-1}(\varphi^{t-1},y_{t-1})\\
=&\inf_{(\psi_{t-1},\psi^{t})\in\boldsymbol\Psi^{t-1}}\int_{E_Y}\cdots\int_{E_Y}e^{\gamma \sum_{k=t-1}^{T-1}\beta^k r_k(x_k,\varphi_k(y_k))+ \gamma\beta^Tr_T(x_T)} \\
& \qquad \qquad \qquad \qquad \prod_{k=t}^T Q(\dif y_k|k-1,y_{k-1},\varphi_{k-1}(y_{k-1}),\psi_{k-1}(y_{k-1}))\\
\geq&\inf_{(\psi_{t-1},\psi^{t})\in\boldsymbol\Psi^{t-1}}\int_{E_Y}e^{\gamma \beta^{t-1} r_{t-1}(x_{t-1},\varphi_{t-1}(y_{t-1}))}V_{t}(y_t)
Q(\dif y_t|t-1,y_{t-1},\varphi_{t-1}(y_{t-1}),\psi_{t-1}(y_{t-1}))\\
=&\inf_{\theta\in\tau(t-1,c)}\int_{E_Y}e^{\gamma \beta^{t-1} r_{t-1}(x_{t-1},\varphi_{t-1}(y_{t-1}))}V_{t}(y_t)
Q(\dif y_t|t-1,y_{t-1},\varphi_{t-1}(y_{t-1}),\psi_{t-1}(y_{t-1}))\\
=&V_{t-1}(y_{t-1}).
\end{align*}
Next, fix $\epsilon>0$, and let $\psi^{t,\epsilon}$ denote an $\epsilon$-optimal selectors sequence starting at time $t$, namely
\begin{align*}
\int_{E_Y}\cdots\int_{E_Y}e^{\gamma \sum_{k=t}^{T}\beta^k r_k(x_k,\varphi_k(y_k))} &\prod_{k=t+1}^T
Q(\dif y_k|k-1,y_{k-1},\varphi_{k-1}(y_{k-1}),\psi^{t,\epsilon}_{k-1}(y_{k-1}))\\
\leq&U_{t}(\varphi^{t},y_{t})+\epsilon.
\end{align*}
Consequently, for any $y_{t-1}\in E_Y$,
\begin{align*}
U_{t-1}&(\varphi^{t-1},y_{t-1})=\inf_{(\psi_{t-1},\psi^{t})\in\boldsymbol\Psi^{t-1}}\int_{E_Y}\cdots\int_{E_Y}
e^{\gamma \sum_{k=t-1}^{T}\beta^k r_k(x_k,\varphi_k(y_k))}\\
&\qquad\qquad\prod_{k=t}^TQ(\dif y_k|k-1,y_{k-1},\varphi_{k-1}(y_{k-1}),\psi_{k-1}(y_{k-1}))\\
&\leq\inf_{\psi_{t-1}\in\tau(t-1,c)}\int_{E_Y}\cdots\int_{E_Y}e^{\gamma \sum_{k=t-1}^{T}\beta^k r_k(x_k,\varphi_k(y_k))} \\
&\quad \prod_{k=t+1}^TQ(\dif y_k|k-1,y_{k-1},\varphi_{k-1}(y_{k-1}),\psi^{t,\epsilon}_{k-1}(y_{k-1}))\cdots
 Q(\dif y_t|t-1,y_{t-1},\varphi_{t-1}(y_{t-1}),\psi_{t-1}(y_{t-1}))\\
&\leq\inf_{\varphi_{t-1}\in\tau(t-1,c)}\int_{E_Y}U_{t}(\varphi^{t},y_t)Q(\dif y_t|t-1,y_{t-1},\varphi_{t-1}(y_{t-1}),\psi_{t-1}(y_{t-1}))+\epsilon\\
&=\inf_{\varphi_{t-1}\in\tau(t-1,c)}\int_{E_Y}V_{t}(y_t)Q(\dif y_t|t-1,y_{t-1},\varphi_{t-1}(y_{t-1}),\psi_{t-1}(y_{t-1}))+\epsilon\\
& =V_{t-1}(y_{t-1})+\epsilon.
\end{align*}
Since $\epsilon$ is arbitrary, \eqref{V} is justified. In particular, $U_{t}(\varphi^{t},\cdot)$ is l.s.a. for any $t\in\cT'$. Finally, in view of \eqref{eq:lemma34-3}, the equality \eqref{eq:lemma34} follows immediately.
This concludes the proof.
\end{proof}

Now we are in the position to prove the main result in this paper.
\begin{theorem}\label{th:main} For  $t=0,\ldots,T$, we have that
\begin{equation}\label{opt}
U^*_{t}\equiv W_{t}.
\end{equation}
Moreover, the policy $\varphi^*$ derived in Lemma \ref{lemma:WTusc} is adaptive robust-optimal, that is
\begin{equation}\label{opt1}
U^*_{t}(y)=U_{t}(\varphi^{t,*},y), \quad t=0,\ldots,T-1.
\end{equation}
\end{theorem}

\begin{proof}
We proceed similarly as in the proof of \cite[Theorem 2.1]{Iyengar2005}, and via backward induction in $t=T,T-1,\ldots,1,0$.

For $t=T$, clearly, $U^*_{T}(y)=W_{T}(y)=e^{\gamma \beta^T r_T(x)}$ for all $y\in E_Y$. For $t=T-1$ we have, for $y\in E_Y$,
\begin{align*}
U^*_{T-1}(y)
& = \sup_{\varphi^{T-1}=\varphi_{T-1}\in {\mathcal A}^{T-1}}\inf_{ \theta \in \tau(T-1,c)} \int_{E_Y}e^{\gamma \beta^{T-1} r_{T-1}(x,\varphi_{T-1}(y))}  W_{T}(y') \,Q( \dif y'\mid T-1,y_{T-1},\varphi_{T-1}(y),\theta )\\
&= \max_{a\in {A}}\inf_{ \theta \in \tau(T-1,c)}\int_{E_Y} e^{\gamma \beta^{T-1} r_{T-1}(x,a)}  W_{T}(y') \, Q( \dif y'\mid T-1,y,a,\theta )\\
& =W_{T-1}(y).
\end{align*}
From the above, using Lemma \ref{lemma:WTusc}, we obtain that $U^*_{T-1}$ is l.s.c. and bounded.

For $t=T-2,\ldots,1,0$, assume that $U^*_{t+1}$ is l.s.c. and bounded. Recalling the notation $\varphi^t=(\varphi_t,\varphi^{t+1})$, we thus have, $y\in E_Y$,
\begin{align*}
U^*_{t}(y)
&=\sup_{(\varphi_t,\varphi^{t+1})\in {\mathcal A}^{t}}\inf_{ \theta \in \tau(t,c)}
	\int_{E_Y}  e^{\gamma\beta^t r_t(x,\varphi_t(y))} U_{t+1}(\varphi^{t+1},y') \,  Q( \dif y'\mid t,y,\varphi_t(y),\theta )\\
&\leq \sup_{(\varphi_t,\varphi^{t+1})\in {\mathcal A}^t}\inf_{ \theta \in \tau(c_t,t)}
\int_{E_Y} e^{\gamma\beta^t r_t(x,\varphi_t(y))} U^*_{t+1}(y')  \,Q( \dif y'\mid t,y,\varphi_t(y),\theta )\\
& = \max_{a\in {A}}\inf_{ \theta \in \tau(t,c)}
	\int_{E_Y}  e^{\gamma\beta^t r_t(x,a)} U^*_{t+1}(y') \, Q( \dif y\mid t,y_t,a,\theta )\\
& = \max_{a\in {A}}\inf_{ \theta \in \tau(t,c)}\int_{E_Y} e^{\gamma\beta^t r_t(y,a)} W_{t+1}(y') \, Q( \dif y'\mid t,y,a,\theta ) \\
& =W_{t}(y).
\end{align*}
Now, fix $\epsilon >0$, and let $\varphi^{t+1,\epsilon}$ denote an $\epsilon$-optimal control strategy starting at time $t+1$, that is
\[
	U_{t+1}(\varphi^{t+1,\epsilon},y)\geq U^*_{t+1}(y)-\epsilon,\quad y\in E_y.
\]
Then,   for $y\in E_Y$, we have
\begin{align*}
U^*_{t}(y)
& = \sup_{(\varphi_t,\varphi^{t+1})\in {\mathcal A}^{t}}\inf_{ \theta \in \tau(t,c)}
		\int_{E_Y}  e^{\gamma\beta^t r_t(x,\varphi_t(y))} U_{t+1}(\varphi^{t+1},y') \, Q( \dif y'\mid t,y,\varphi_t(y),\theta )\\
&\geq \sup_{(\varphi_t,\varphi^{t+1})\in {\mathcal A}^{t}}\inf_{ \theta \in \tau(t,c)}
		\int_{E_Y}  e^{\gamma\beta^t r_t(x,\varphi_t(y))}  U_{t+1}(\varphi^{t+1,\epsilon},y') \, Q( \dif y'\mid t,y,\varphi_t(y),\theta )\\
& \geq   \max_{a\in A}\inf_{ \theta \in \tau(t,c)}\int_{E_Y} 	 e^{\gamma\beta^t r_t(x,a)} U^*_{t+1}(y') 	\, Q( \dif y'\mid t,y,a,\theta) -\epsilon  \\
& =   \max_{a\in A}\inf_{ \theta \in \tau(t,c)}\int_{E_Y}  W_{t+1}(y')\, 	Q(\dif y'\mid t,y,a, \theta )-\epsilon\\
&=W_{t}(y)-\epsilon.
\end{align*}
Since $\epsilon$ was arbitrary, the proof of \eqref{opt} is done. In particular, we have that for any $t\in\cT$, the function $U^*_{t}(\cdot)$ is l.s.c. as well as bounded.

It remains to justify the validity of equality \eqref{opt1}. We will proceed again by (backward) induction in $t$.
For $t=T-1$, using \eqref{opt}, we have that
\begin{align*}
U^*_{T-1}(y) &
= W_{T-1}(y) =e^{\gamma \beta^{T-1}r_{T-1}(x,\varphi^*_{T-1}(y))}
\inf_{ \theta \in \tau(t,c)}\int_{E_Y} e^{\gamma \beta ^Tr_T(x')}\, Q(\dif y'\mid T-1,y,\varphi^*_{T-1}(y),\theta ) \\
& =e^{\gamma \beta^{T-1}r_{T-1}(x,\varphi^*_{T-1}(y))}\inf_{\mathbb{Q}\in {\mathcal Q}^{\varphi^{T-1,*}, {\boldsymbol \Psi}^{T-1}}_{y,T-1}} \left (\bE_{\mathbb Q} e^{\gamma \beta^Tr_T(X_T)}\right )\\
& =U_{T-1}(\varphi^{T-1,*},y).
\end{align*}
Moreover, by Lemma~\ref{lemma:U-lsa}, we get that
$$
U^*_{T-1}(y)=U_{T-1}(\varphi^{T-1,*},y)=\bE_{\mathbb Q^{\varphi^{T-1,*},\psi^{T-1,*}}_{y, T-1}}e^{\gamma \beta^{T-1}r_{T-1}(x,\varphi^*_{T-1}(y)) + \gamma \beta^Tr_T(X_T)}.
$$
For $t=T-2$, using again \eqref{opt}, Lemma~\ref{lemma:WTusc}, and Lemma~\ref{lemma:U-lsa}, we have
\begin{align*}
U^*_{T-2}(y)
& = W_{T-2}(y)=e^{\gamma \beta^{T-2}r_{T-2}(x,\varphi^*_{T-2}(y))}\int_{E_Y} W_{T-1}(y')\,
		Q(\dif y'\mid T-2,y,\varphi^*_{T-2}(y),\psi^*_{T-2}(y,\varphi^*_{T-2}(y)) ) \\ 		
& =e^{\gamma \beta^{T-2}r_{T-2}(x,\varphi^*_{T-2}(y))}\int_{E_Y} U_{T-1}(\varphi^{T-1,*},y')\,
		Q(\dif y'\mid T-2,y,\varphi^*_{T-2}(y),\psi^*_{T-2}(y,\varphi^*_{T-2}(y)) )\\
& = e^{\gamma \beta^{T-2}r_{T-2}(x,\varphi^*_{T-2}(y))}\times \int_{E_Y}\left (\bE_{\mathbb Q^{\varphi^{T-1,*},\psi^{T-1,*}}_{y', T-1}} e^{\gamma \beta^{T-1} r_{T-1}(x'),\varphi^*_{T-1}(y'))+\gamma \beta^Tr_T(X_T)}\right ) \\
& \qquad\qquad\qquad Q(\dif y'\mid T-2,y,\varphi^*_{T-2}(y), \psi^*_{T-2}(y,\varphi^*_{T-2}(y))) \\
& =\bE_{\bQ^{\varphi^{T-2,*},\psi^{T-2,*}}_{y,T-2}}e^{\gamma \beta^{T-2}r_{T-2}(x,\varphi^*_{T-2}(y)) + \gamma \beta^{T-1} r_{T-1}(x'),\varphi^*_{T-1}(y'))+\gamma \beta^Tr_T(X_T)}.
\end{align*}
Hence, we have that $U^*_{T-2}(y)$ is attained at $\varphi^{T-2,*}$, and therefore $U^*_{T-2}(y)=U_{T-2}(\varphi^{T-2,*},y)$.
The rest of the proof of \eqref{opt1} proceeds in an analogous way. The proof is complete.
\end{proof}

\section{Exponential Discounted Tamed Quadratic Criterion Example}\label{sec:ex}

In this section, we consider a linear quadratic control problem under model uncertainty as a numerical demonstration of the adaptive robust method.
To this end, we consider the 2-dimensional controlled process
\begin{align*}
X_{t+1} = B_1X_t+B_2\varphi_t+Z_{t+1},
\end{align*}
where $B_1$ and $B_2$ are two $2\times2$ matrices and $Z_{t+1}$ is a 2-dimensional normal random variable with mean 0 and convariance matrix
$$
\Sigma^*=\begin{pmatrix}
\sigma^{*,2}_1 & \sigma^{*,2}_{12} \\
\sigma^{*,2}_{12} & \sigma^{*,2}_2
\end{pmatrix},
$$
where $\sigma^{*,2}_1$, $\sigma^{*,2}_{12}$, and $\sigma^{*,2}_2$ are unknown.
Given observations $Z_1,\ldots,Z_t$, {we consider an} unbiased estimator, say $\widehat\Sigma_t=\begin{pmatrix}
\widehat{\sigma}^2_{1,t} & \widehat{\sigma}^2_{12,t} \\
\widehat{\sigma}^2_{12,t} & \widehat{\sigma}^2_{2,t}
\end{pmatrix}$, of the covariance matrix $\Sigma^*$, given as
\begin{align*}
\widehat\Sigma_t=\frac{1}{t+1}\sum_{i=1}^tZ_iZ^\top_i,
\end{align*}
which can be updated recursively as
\begin{align*}
\widehat\Sigma_t=\frac{t(t+1)\widehat\Sigma_{t-1}+tZ_tZ_t^\top}{(t+1)^2}.
\end{align*}
With slight abuse of notations, we denote by $\Sigma$, $\Sigma^*$, and $\widehat\Sigma_t$ the column vectors
\begin{align*}
\Sigma^\top &=(\sigma^2_1,\sigma^2_{12},\sigma^2_2)\\
\Sigma^{*,\top} &=(\sigma^{*,2}_1,\sigma^{*,2}_{12},\sigma^{*,2}_2)\\
\widehat\Sigma^\top_t&=(\widehat{\sigma}^2_{1,t},\widehat{\sigma}^2_{12,t},\widehat{\sigma}^2_{2,t}).
\end{align*}
The corresponding parameter set is defined as
$$
\boldsymbol\Theta:=\left\{\Sigma^\top=(\Sigma_1,\Sigma_{12},\Sigma_2)\in\bR^3:\ 0\leq\Sigma_1,\ \Sigma_2\leq\overline{\Sigma},\ \Sigma_{12}^2\leq\Sigma_1\Sigma_2\right\},
$$
where $\overline{\Sigma}$ is some {fixed} positive constant. Note that the set $\boldsymbol\Theta$ is a compact subset of $\bR^3$.

Putting the above together and considering the augmented state process $Y_t=(X_t,\widehat\Sigma_t)$, $t\in\cT$, and some finite control set $A\subset\bR^2$, we get that the function $S$ defined in \eqref{eq:mm} is given by
$$
S(x,a,z)=B_1x+B_2a+z, \quad x,z\in\bR^2,\ a\in A,
$$
and the function $R(t,c,z)$ showing in \eqref{eq:R} satisfies that
\begin{align*}
R(t,c,z)=(\bar c_1,\bar c_2,\bar c_3)^\top, \quad
\begin{pmatrix}
\bar{c}_1 & \bar{c}_3 \\
\bar{c}_3 & \bar{c}_2
\end{pmatrix}=
\frac{(t+1)(t+2)
\begin{pmatrix}
c_1 & c_3 \\
c_3 & c_2
\end{pmatrix}
+(t+1)zz^\top
}{(t+2)^2},
\end{align*}
where $z\in\bR^2$, $t\in\cT'$, $c=(c_1,c_2,c_3)$.
Then, function $\mathbf{G}$ defined in \eqref{eq:T} is specified accordingly.

It is well-known that $\sqrt{t+1}(\widehat\Sigma_t-\Sigma^*)$ converges weakly to 0-mean normal dsitribution with covariance matrix
$$
M_{\Sigma}=\begin{pmatrix}
2\sigma^{*,4}_1 & 2\sigma^{*,2}_1\sigma^{*,2}_{12} & 2\sigma^{*,4}_{12} \\
2\sigma^{*,2}_1\sigma^{*,2}_{12} & \sigma^{*,2}_1\sigma^{*,2}_2+\sigma^{*,4}_{12} & 2\sigma^{*,2}_{12}\sigma^{*,2}_2 \\
2\sigma^{*,4}_{12} & 2\sigma^{*,2}_{12}\sigma^{*,2}_2 & 2\sigma^{*,4}_2
\end{pmatrix}.
$$
We replace every entry in $M_{\Sigma}$ with the corresponding estimator at time $t\in\cT'$ and denote by $\widehat{M}_t(\widehat\Sigma_t)$ the resulting matrix.
With probability one, the matrix $\widehat{M}_t(\widehat\Sigma_t)$ is positive-definite.
Therefore, we get the confidence region for $\sigma^{*,2}_1$, $\sigma^{*,2}_{12}$, and $\sigma^{*,2}_2$ as
\begin{align*}
\tau(t,c)=\left\{\Sigma\in\boldsymbol\Theta:(t+1)(\Sigma-c)^\top \widehat{M}^{-1}_t(c)(\Sigma-c)\leq\kappa\right\},
\end{align*}
where $\kappa$ is the $1-\alpha$ quantile of $\chi^2$ distribution with 3 degrees of freedom for some confidence level $0<\alpha<1$.

We further take functions $r_T(x)=\min\{b_1,\max\{b_2,x^\top K_1x\}\}$ and $r_t(x,a)=\min\{b_1,\max\{b_2,x^\top K_1x+a^\top K_2a\}\}$, $t\in\cT'$, where $x,a\in\bR^2$, $b_1>0$, $b_2<0$, and $K_1$ and $K_2$ are two fixed 2-by-2 matrices with negative trace.


{For this example, all conditions of the adaptive robust framework of Section~\ref{sec:robust} are easy to verify, except for the u.h.c. property of set-valued function $\tau(t,\cdot)$, which we establish in the following lemma}.

\begin{lemma}\label{lemma:hemi}
	For any $t\in\cT'$, the set valued function $\tau(t,\cdot)$ is upper hemi-continuous.
\end{lemma}

\begin{proof}
	Fix any $t\in\cT'$ and $c_0\in\boldsymbol\Theta$.
	According to our earlier discussion, the matrix $\widehat{M}_t(c_0)$ is positive-definite. Hence, its inverse admits the Cholesky decomposition $\widehat{M}^{-1}_t(c_0)=L_t(c_0)L^{\top}_t(c_0)$. Consider the change of coordinate system via the linear transformation $\cL c=L^\top_t(c_0)c$, and we name it system-$\cL$. Let $E\subset\boldsymbol\Theta$ be open and such that $\tau(t,c_0)\subset E$. Note that $\cL\tau(t,c_0)$ is a closed ball centered at $\cL c_0$ in the system-$\cL$. Also, the mapping $\cL$ is continuous and one-to-one, hence $\cL E$ is an open set and $\cL\tau(t,c_0)\subset\cL E$. Then, we have that there exists an open ball $B_r(\cL c_0)$ in the system-$\cL$ centered at $\cL c_0$ with radius $r$ such that $\cL\tau(t,c_0)\subset B_r(\cL c_0)\subset\cL E$.
	
	Any ellipsoid centered at $c'$ in the original coordinate system has representation $(c-c')^\top F(c-c')=1$ which can be written as $(L^\top_tc-L^\top_tc')L^{-1}F(L^\top)^{-1}(L^\top c-L^\top c')=1$.
	Hence, it is still an ellipsoid in the $\cL$-system after transformation. To this end, we define on $\boldsymbol\Theta$ a function $h(c):=\|\cL c-\cL c_0\|+\max\{r_i(c),i=1,2,3\}$, where $\|\cdot\|$ is the Euclidean norm in the system-$\cL$, and $r_i(c)$, $i=1, 2, 3$, are the lengths of the three semi axes of the ellipsoid $\cL \tau(t,c)$. It is clear that $r_i(c)$, $i=1,2,3$ are continuous functions.
	
	Next, it is straightforward to check that $f$ is a non-constant continuous function. Therefore, we consider the set $D:=\{c\in\boldsymbol\Theta: h(c)<r\}$ and see that it is an open set in $\boldsymbol\Theta$ and non-empty as $c_0\in D$.
	Moreover, for any $c\in D$, we get that the ellipsoid $\cL\tau(t,c)\subset B_r(\cL c_0)$.
	Hence, $\tau(t,c)\subset E$, and we conclude that $\tau(t,\cdot)$ is u.h.c..
\end{proof}

Thus, according to Theorem~\ref{th:main}, the dynamic risk sensitive optimization problem under model uncertainty can be reduced to the Bellman equations given in \eqref{eq:bellmanEquationrobustIII}:
\begin{align}
  W_{T}(y) & = e^{\gamma \beta^T r_T(x)}, \label{eq:exmp1} \\
W_{t}(y) & = \sup_{a\in A} \inf_{ \theta \in \tau(t,c)}
\int_{\bR^2}W_{{t+1}}(\mathbf{G}(t,y,a,z) )e^{\gamma \beta^t (r_t(x,a))}
   f_Z(z;\theta)dz,\label{eq:exmp2}\\
   y&=(x,c_1,c_2,c_3)\in E_Y, \   t=T-1, \ldots, 0, \nonumber
\end{align}
where $f_Z(\cdot;\theta)$ is the density function for two dimensional normal random variable with mean 0 and covariance parameter $\theta$.
In the next section, {using \eqref{eq:exmp1}-\eqref{eq:exmp2}}, we will  compute numerically $W_t$ by  a machine learning based method. Note that the dimension of the state space $E_Y$ is five in the present case, for which the traditional grid-based numerical method becomes extremely inefficient. Hence, we employ the new approach introduced in \cite{ChenLudkovski2019} to overcome the challenges met in our high dimensional robust stochastic control problem.

\section{Machine Learning Algorithm and Numerical Results}\label{sec:ml}

In this section, we describe our machine learning based method and present the numerical results for our example.
Similarly to \cite{ChenLudkovski2019}, we discretize the state space the relevant state space in the spirit of the regression Monte Carlo method and adaptive design by creating a random (non-gridded) mesh for the process $Y=(X,C)$. Note that the component $X$ depends on the control process, hence at each time $t$ we randomly select from the set $A$ a value of $\varphi_t$, and we randomly generate a value of $Z_{t+1}$, so to simulate the value of $X_{t+1}$. Next, for each $t$, we construct the convex hull of simulated $Y_t$ and uniformly generate in-sample points from the convex hull to obtain a random mesh of $Y_t$. Then, we solve the equations \eqref{eq:exmp1}--\eqref{eq:exmp2}, and compute the optimal trading strategies at all mesh points.

The key idea of our machine learning based method is to utilize a non-parametric value function approximation strategy called Gaussian process surrogate. For the purpose of solving the Bellman equations \eqref{eq:exmp1}--\eqref{eq:exmp2}, we build GP regression model for the value function $W_{t+1}(\cdot)$ so that we can evaluate
$$
\int_{\bR^2}W_{{t+1}}(\mathbf{G}(t,y,a,z) )e^{\gamma \alpha^t (r_t(x,a))}
   f_Z(z;\theta)dz.
$$
We also construct GP regression model for the optimal control $\varphi^*$. It permits us to apply the optimal strategy to out-of-sample paths without actual optimization, which allows for a significant reduction of the computational cost.

As the GP surrogate for the value function $W_t$ we consider a regression model $\widetilde{W}_t(y)$ such that for any $y^1,\ \ldots,\ y^N\in E_Y$, with $y^i\neq y^j$ for $i\neq j$, the random variables $\widetilde{W}_t(y^1),\ \ldots,\ \widetilde{W}_t(y^N)$ are jointly normally distributed. Then, given training data $(y^i,W_t(Y^i))$, $i =1,\ \ldots,\ N$, for any $y\in E_Y$, the predicted value $\widetilde{W}_t(y)$, providing an estimate (approximation) of $W_t(y)$ is given by
\begin{align*}
  \widetilde{W}(y) = \left(k(y,y^1),\ldots,k(y,y^N)\right)[\mathbf{K}+\epsilon^2\mathbf{I}]^{-1}\left(W_t(y^1),\ldots,W_t(y^N)\right)^T,
\end{align*}
where $\epsilon$ is a tuning parameter, $\mathbf{I}$ is the $N\times N$ identity matrix and the matrix $\mathbf{K}$ is defined as $\mathbf{K}_{i,j}=k(y^i,y^j)$, $i,\ j=1,\ \ldots,\ N$. The function $k$ is the kernel function for the GP model, and in this work we choose the kernel as the Matern-5/2. Fitting the GP surrogate $\widetilde{W}_t$ means to estimate the hyperparameters inside $k$ through the training data $(y^i,W_t(y^i))$, $i=1,\ \ldots,\ N$ for which we take $\epsilon=10^{-5}$. The GP surrogates for $\varphi^*$ is obtained in an analogous way.

Given the mesh points $\{y^i_t,\ i=1,\ \ldots,\ N_t,\ t\in\cT'\}$, the overall algorithm proceeds as follows:\\
\textit{Part A:} Time backward recursion for $t=T-1,\ldots, 0$.
\begin{enumerate}
  \item Assume that $W_{t+1}(y^i_{t+1})$, and $\varphi^*_{t+1}(y^i_{t+1})=(\varphi^{1,*}_{t+1}(y^i_{t+1}),\varphi^{2,*}_{t+1}(y^i_{t+1}))$, $i=1,\ldots,N_t$, are numerically approximated as $\overline W_{t+1}(y^i_{t+1})$, $\overline{\varphi}^{1,*}_{t+1}(y^i_{t+1})$ and $\overline{\varphi}^{2,*}_{t+1}(y^i_{t+1})$, $i=1,\ldots,N_t$, respectively. Also suppose that the corresponding GP surrogates $\widetilde{W}_{t+1}$, $\widetilde{\varphi}^{1,*}_{t+1}$, and $\widetilde{\varphi}^{2,*}_{t+1}$ are fitted through training data $(y^i_{t+1},\overline{W}_{t+1}(y^i_{t+1}))$, $(y^i_{t+1},\overline{\varphi}^{1,*}_{t+1}(y^i_{t+1}))$, and $(y^i_{t+1},\overline{\varphi}^{2,*}_{t+1}(y^i_{t+1}))$, $i=1,\ldots,N_t$, respectively.
  \item For time $t$, any $a\in A$, $\theta\in\tau(t,c)$ and each $y^i_t$, $i=1,\ \ldots, N_t$, use one-step Monte Carlo simulation to estimate the integral
 \begin{align*}
	w_t(y,a,\theta)=\int_{\bR^2}W_{{t+1}}(\mathbf{G}(t,y,a,z) )e^{\gamma \alpha^t (r_t(x,a))}
   f_Z(z;\theta)dz.
\end{align*}
For that, if $Z^1_{t+1},\ \ldots,\ Z^M_{t+1}$ is a sample of $Z_{t+1}$ drawn from the normal distribution corresponding to parameter $\theta$, where $M>0$ is a positive integer, then estimate the above integral as
\begin{align*}
	\widetilde{w}_t(y,a,\theta)=\frac{1}{M}\sum_{i=1}^M\widetilde{W}_{{t+1}}(\mathbf{G}(t,y,a,Z^i_{t+1}) )e^{\gamma \alpha^t (r_t(x,a))}.
\end{align*}
	\item For each $y^i_t$, $i=1,\ \ldots,\ N_t$, and any $a\in A$, compute
$$
\overline w_{t}(y^i_t,a)=\inf_{\theta\in\tau(t,c)}\widetilde{w}_{t}(y^i_t,a,\theta).
$$

\item Compute
$$
\overline W_t(y_t^i)=\max_{a\in A}\overline{w}_t(y^i_t,a),
$$
and obtain a maximizer $\overline{\varphi}^*_t(y_t^i)=(\overline{\varphi}^{1,*}_t(y_t^i),\overline{\varphi}^{2,*}_t(y_t^i))$, $i=1,\ldots,N_t$.
\item Fit a GP regression model for $V_t(\,\cdot\,)$ using the results from Step~4 above. Fit GP models for $\varphi^{1,*}_t(\,\cdot\,)$ and $\varphi^{2,*}_t(\,\cdot\,)$ as well; these are needed for obtaining values of the optimal strategies for out-of-sample paths in Part B of the algorithm.
	\item Goto 1: Start the next recursion for $t-1$.
\end{enumerate}
\textit{Part B:} Forward simulation to evaluate  the performance of the GP surrogates $\varphi^{1,*}_t(\,\cdot\,)$ and $\varphi^{2,*}_t(\,\cdot\,)$, $t=0, \ldots,T-1$, over the out-of-sample paths.
\begin{enumerate}
  \item Draw $K>0$ samples of i.i.d.  $Z_{1}^{*,i},\ldots,Z_{T}^{*,i}$, $i=1,\ldots,K$, from the normal distribution corresponding to the assumed true parameter $\theta^*$.
  \item All paths will start from the initial state $y_0$. The state along each path $i$ is updated according to $\mathbf G(t,y_t^i,\widetilde{\varphi}^*_t(y_t^i),Z_{t+1}^{*,i})$, where $\widetilde{\varphi}^*_t=(\widetilde\varphi^{1,*}_t,\widetilde\varphi^{2,*}_t)$ is the GP surrogate fitted in Part A. Also, compute the running reward $r_t(x^i_t,\widetilde{\varphi}^*_t(y^i_t))$.
  \item Obtain the terminal reward $r_T(x^i_T)$, generated by $\widetilde{\varphi}^*$ along the path corresponding to the sample of  $Z_{1}^{*,i},\ \ldots,\ Z_{T}^{*,i}$, $i=1,\ \ldots,\ K$, and compute
           \begin{align}\label{eq:ar}
       W^{\text{ar}}:=\frac{1}{\gamma}\ln\left(\frac{1}{K}\sum_{i=1}^{K}e^{\gamma(\sum_{t=0}^{T-1}\beta^tr_t(x^i_t,\widetilde{\varphi}^*_t(y^i_t))+\beta^Tr_T(x^i_T))}\right)
      \end{align}
      as an estimate of the performance of the optimal adaptive robust risk sensitive strategy $\varphi^*$.
\end{enumerate}

For comparison, we also analyze the optimial risk sensitive strategies of the adaptive and strong robust control methods. In \eqref{eq:exmp2}, if we take $\tau(t,c)=\{c\}$ for any $t$, then we obtain the adaptive risk sensitive strategy. On the other hand, by taking $\tau(t,c)=\boldsymbol{\Theta}$ for any $t$ and $c$, we get the strong robust strategy.
We will compute $W^{\text{ad}}$ and $W^{\text{sr}}$ the risk sensitive criteria of adaptive and strong robust, respectively, in analogy to \eqref{eq:ar}.

Next, we apply the machine learning algorithm described above by solving \eqref{eq:exmp1}--\eqref{eq:exmp2} for a specific set of parameters.
In particular, we take: $T=10$ with one period of time corresponding to one-tenth of a year;
the discount factor being equal to 0.3 or equivalently $\beta=0.3$; the initial state $X^\top_0=(2,2)$; the confidence level $\alpha=0.1$; in Part A of our algorithm the number of one-step Monte Carlo simulations is $M=100$; the number of forward simluations in Part B is taken $K=2000$; the control set $A$ is approximated by the compact set $[-1,1]^2$; the relevant matrices are
$$
B_1=B_2=\begin{pmatrix}
0.5 & -0.1\\
-0.1 & 0.5
\end{pmatrix},
\quad
K_1=\begin{pmatrix}
0.7 & -0.2 \\
-0.2 & 0.7
\end{pmatrix},
\quad
K_2=\begin{pmatrix}
-200 & 100 \\
100 & -200
\end{pmatrix}.
$$
The assumed true covariance matrix for $Z_t$, $t\in\cT$, as well as initital guess are
$$
\Sigma^*=\begin{pmatrix}
0.009 & 0.006 \\
0.006 & 0.016
\end{pmatrix},
\quad
\widehat\Sigma_0=\begin{pmatrix}
0.00625 & 0.004 \\
0.004 & 0.02025
\end{pmatrix},
$$
respectively.
The parameter set is chosen as $\boldsymbol\Theta=\tau(0,c_0)$, where $c_0^\top=(0.00625,0.004,0.02025)$. For all three control approaches, we compute $W^{\text{ar}}$, $W^{\text{ad}}$, and $W^{\text{sr}}$, respectively, for the risk sensitive parameters $\gamma=0.2$ and $\gamma=1.5$.

Finally, we report on the computed values of the optimality criterion corresponding to  three different methods: adaptive robust (AR), adaptive (AD) and strong robust (SR).

\begin{table}[h!]
\centering
\renewcommand{\arraystretch}{1.3}
\begin{tabular}{c|ccc}
\hline
  \multicolumn{1}{c}{ }& \multicolumn{1}{|c}{$W^{\textrm{ar}}$} & $W^{\textrm{ad}}$ & $W^{\textrm{sr}}$ \\
  \hline
   $\gamma=0.2$ & \multicolumn{1}{|c}{-319.81} & -323.19 & -329.53\\
 $\gamma=1.5$ & \multicolumn{1}{|c}{-427.76} & -427.97 & -442.97 \\
 \hline
\end{tabular}
\bigskip
\caption{Risk sensitive criteria for AR, AD, and SR.}
\label{table:comparison}
\end{table}

%
%

\bibliographystyle{alpha}
\bibliography{MathFinanceMaster-01-18-2021}

\end{document}